\documentclass[12pt]{article}
\usepackage{mathrsfs}
\usepackage{amsmath}
\usepackage{amssymb}
\usepackage{amsthm}
\usepackage{amsfonts}
\usepackage{amscd}
\usepackage[mathscr]{eucal}
\newcommand{\Z} {{\mathbb  Z}}
\newcommand{\Q}{{\mathbb  Q}}
\newcommand{\F}{{\mathbb  F}}
\newcommand{\C}{{\mathbb  C}}

\newcommand{\R} {{\mathbb R}}

\begin{document}
\parindent  25pt
\baselineskip  10mm
\textwidth  15cm    \textheight  23cm \evensidemargin -0.06cm
\oddsidemargin -0.01cm

\title{ {On quadratic twists of elliptic curves and some
applications of a refined version of Yu's formula
 }}
\author{\mbox{}
{ Derong Qiu }
\thanks{ \quad E-mail:
derong@mail.cnu.edu.cn, \ derongqiu@gmail.com } \\
(School of Mathematical Sciences,
 Capital Normal University, \\
Beijing 100048, P.R.China )  }

\date{}
\maketitle
\parindent  24pt
\baselineskip  10mm
\parskip  0pt

\par   \vskip 0.4cm

{\bf Abstract} \quad In this paper, we study some cohomology groups
and quadratic twists of elliptic curves, and apply Tate local duality and
the results of Kramer-Tunnell on local norm cokernel to give
a refined version of Yu's formula in the case of
elliptic curves. Then, by using this refinement formula, we obtain explicit
orders of Shafarevich-Tate groups of some elliptic curves in quadratic number
fields, including a few unconditional cases.
\par  \vskip  0.2 cm

\par  \vskip  0.1 cm


\par     \vskip  0.4 cm

\hspace{-0.8cm}{\bf 0. Introduction}

\par \vskip 0.2 cm

Let $ K/F $ be a finite Galois extension of number fields with Galois group $ G, \
A $ an abelian variety defined over $ F, $ and let $ \amalg\hskip-7pt\amalg(A/F),
 \amalg\hskip-7pt\amalg(A/K) $ denote the Shafarevich-Tate groups of $ A $ over $ F $
and $ K, $ respectively. In [Y], under the assumption that these groups are finite,
Yu computed the value
$ \sharp \amalg\hskip-7pt\amalg(A/K)^{G} / \sharp \amalg\hskip-7pt\amalg(A/F). $
In particular, when $ K $ is a quadratic extension of $ F, $ Yu derived a formula
relating $ \amalg\hskip-7pt\amalg(A/F), \amalg\hskip-7pt\amalg(A/K) $ and
$ \amalg\hskip-7pt\amalg(A^{\chi }/F) $ as follows, where $ A^{\chi } $ is the twist of
$ A $ by the non-trivial character $ \chi $ of $ G. $
\par \vskip 0.2 cm

{\bf Yu's formula} (see [Y], Main Theorem, p.212) \ Assume
that the Shafarevich-Tate groups are finite. Let $ A^{\prime } $ be the dual abelian
variety of $ A. $ Then
$$ \frac{\sharp \amalg\hskip-6pt\amalg(A/F) \cdot \sharp
\amalg\hskip-6pt\amalg(A^{\chi }/F)}{\sharp \amalg\hskip-6pt\amalg(A/K)}
= \frac{\sharp \widehat{H}^{0}(G, A^{\prime }(K)) \cdot
\sharp H^{1}(G, A(K))}{\prod _{v \in M_{F}} \sharp H^{1}(G_{v_{K}}, A(K_{v_{K}}))}, $$
where $ v_{K} $ is the fixed place of $ K $ lying above $ v $ for each $ v \in M_{F}. $

A further question is to determine these quantities of the right-hand side of
this equality, which seems not too easy to practically calculate in general.

In this paper, for the case of elliptic curves, we refine the above Yu's formula
in an elementary and simpler form. In fact, by studying some cohomology groups
of quadratic twists of elliptic curves, we can relate the order of
$ H^{1}(G, A(K)) $ with the groups of rational points of $ A $ and its twist
(see Theorem 1.5 in the following Section 1), then by applying Tate local duality
(see [Ta1]), and the results
of Kramer and Tunnell on local norm cokernel (see Prop.$1\sim$5 in [Kr] and Thm.7.6 in[KT]),
we can in some extent calculate the order of
local cohomology groups $ H^{1}(G_{v_{K}}, A(K_{v_{K}})). $ Our refined version of Yu's
formula in the case of elliptic curves is as follows

{ \bf A refined version of Yu's formula } (see Theorem 2.1 in the following). \ Assume
that the Shafarevich-Tate groups $ \amalg\hskip-6pt\amalg(E/F),
\amalg\hskip-6pt\amalg(E_{D}/F) $ and $ \amalg\hskip-6pt\amalg(E/K) $
are finite. Then
$$ \frac{\sharp \amalg\hskip-6pt\amalg(A/F) \cdot \sharp
\amalg\hskip-6pt\amalg(A_{D}/F)}{\sharp \amalg\hskip-6pt\amalg(A/K)}
= 2^{r_{D, F} - r_{F} - \delta (A, F, K)} \cdot (A(F) :
N(A(K)))^{2}, $$ where the value $ \delta (A, F, K) $ can often be explicitly calculated
(see Theorem 2.1 in the following Section 2 for a precise meaning of
$ \delta (A, F, K)$).

One of our motivation here is, Yu's formula makes it possible to determine
the Shafarevich-Tate group of $ A $ over the quadratic extension $ K, $ when
these groups of $ A $ and its twist over $ F $ are known.
One of the uses of our refined version of this formula is that the
values $ \delta (A, F, K) $ can be usually directly calculated
(e.g., see Lemma 3.2 and the proof of Thm.4.1
in the following). As for applications, by this refined version of Yu's formula, and
using some known results of Shafarevich-Tate groups of several family of elliptic
curves over $ \Q, $ we obtain explicit orders of their Shafarevich-Tate groups
over some quadratic number fields, including a few unconditional
cases (see Theorems 3.3, 4.1 and Corollary 3.4 in the following).
\par     \vskip  0.4 cm

\hspace{-0.8cm}{\bf 1. Quadratic twists and cohomology groups}

\par \vskip 0.2 cm

Let $ F, K $ be number fields with $ K = F ( \sqrt{D}) $ a quadratic
extension of $ F $ for some $ D \in F^{ \ast } \setminus F^{ \ast
^{2}}. $ Let $ G = \text{Gal}(K / F) = < \sigma
> $ be its Galois group with a generator $ \sigma . $ Let
$ E: \ y^{2} = x^{3} + a x + b $ be an elliptic curve defined over $
F. $ Its quadratic $ D-$twist is given by $ E_{D}: \ y^{2} = x^{3} +
a D^{2} x + b D^{3}. $ By the Mordell-Weil Theorem (see
[Si1,Thm.6.7,p.239]), the group $ E(F) $ of $ F-$rational points of $ E $ is a
finitely generated abelian group, so are the groups $ E_{D}(F) $ and
$ E(K). $ For simplicity, in the following, we denote $ r_{F} =
\text{rank} E(F), \ r_{D, F} = \text{rank} E_{D}(F) $ and $ r_{K} =
\text{rank} E(K). $ Let $ \amalg\hskip-7pt\amalg(E/F),
\amalg\hskip-7pt\amalg(E_{D}/F) $ and $ \amalg\hskip-7pt\amalg(E/K)
$ be the Shafarevich-Tate groups of $ E $ over $ F, \ E_{D} $ over $
F $ and $ E $ over $ K $ respectively
(see [Si1, p.332] for the definition). \\
Since $ E $ and $ E_{D} $ are $ K-$isomorphic as given by
$$ \phi _{D}: \ E_{D} \longrightarrow E, \quad (x, y) \longmapsto
\left(\frac{x}{( \sqrt{D})^{2}}, \ \frac{y}{( \sqrt{D})^{3}}
\right),
$$ we have $ E(K) \cong E_{D}(K), \
\amalg\hskip-6pt\amalg(E_{D}/K) \cong \amalg\hskip-6pt\amalg(E/K) $
and $ r_{F} + r_{D, F} = r_{K} $ (see [ABF, p.5], [RS, Lemma 2.1]).
Denote
$$ R_{D}(F) = \phi _{D}(E_{D}(F)) = \left\{ \left(\frac{x}{D}, \
\frac{y}{( \sqrt{D})^{3}} \right): \ (x, y) \in E_{D}(F) \right\}
\cup \{ O \} \subset E(K). $$ Obviously, $ R_{D}(F) $ is a subgroup
of $ E(K), $ and $ R_{D}(F) \cong E_{D}(F) $ as abstract groups.
\par  \vskip 0.2 cm

Throughout this paper, for a set $ S, $ we denote its cardinal by $
\sharp S. $ For arbitrary abelian group $ A $ and positive integer $
m, $ we denote $ m A = \{m a : \ a \in A \} $ and $ A [m] = \{a \in
A : \ m a = 0 \}. $ If $ A $ is a $ G-$module, then one has the
following Tate cohomology groups : \ $ \widehat{\text{H}}^{n}(G, A)
= \text{H}^{n}(G, A) \quad \text{if} \ n \geq 1; \
\widehat{\text{H}}^{0}(G, A) = A^{G} / (1 + \sigma )A. $ For the
basic facts of cohomology groups $ \text{H}^{n} (G, A) \ (0 \leq n
\in \Z)$ and Tate cohomology groups $ \widehat{\text{H}}^{m}(G, A) \
(m \in \Z) $ of $ G-$module $ A, $ see [Se, Chapt.VIII]and [AW, $\S6$].
\par \vskip 0.2 cm

{\bf Lemma 1.1.} \ $ R_{D}(F) = \{ P \in E(K): \ \sigma (P) = - P
\}. $
\par  \vskip 0.2cm

The two maps \ $
\varphi _{1}: \ E (K) \longrightarrow E (K), \quad P
\longmapsto P + \sigma P \quad ( \forall P \in E (K)) $ and $
\varphi _{2}: \ E (K) \longrightarrow E (K), \quad P \longmapsto P
- \sigma P \quad ( \forall P \in E (K)) $
are endomorphisms of abelian group $ E (K) $ with kernels $
\text{ker} \varphi _{1} = R_{D}(F) $ and $
\text{ker} \varphi _{2} = E (F) $ respectively. We denote
$ N(E(K)) = \text{im} \varphi _{1}, $ the images of $ \varphi _{1}; $ and
$ T_{D} (F) = \text{im} \varphi _{2}, $ the images of $ \varphi
_{2}. $ Obviously $ N(E(K)), T_{D}(F) $ and $ R_{D}(F) $ are
finitely generated abelian groups because they are subgroups of $ E
(K). $ We have
$$ 2 E (F) \subset N(E(K)) \subset E (F), \quad 2 R_{D}(F) \subset
T_{D} (F) \subset R_{D}(F). $$
Moreover, by the former discussion and the exact
sequences of abelian groups
$ O \rightarrow R_{D}(F) \rightarrow E (K) \rightarrow N(E(K))
\rightarrow O $ and $
O \rightarrow E(F) \rightarrow E (K) \rightarrow
T_{D} (F) \longrightarrow O, $ we have \
$ \text{rank}E(K) = \text{rank}R_{D}(F) + \text{rank}
N(E(K)) = \text{rank} E(F) + \text{rank} T_{D}(F), \\
\text{rank} T_{D}(F) = \text{rank} E_{D}(F) = \text{rank} R_{D}(F),
\quad \text{rank} N(E(K)) = \text{rank} E(F). $

In particular, the quotient groups $ E(F) / N(E(K)) $ and $
R_{D}(F) / T_{D}(F) $ are finite abelian groups.
\par  \vskip 0.2cm

{\bf Lemma 1.2.} \\
(1) \ $ R_{D}(F)[2] = R_{D}(F) \cap E (F) = E(F)[2]. $ \\
(2) \ $ T_{D}(F)[2] = T_{D}(F) \cap E (F) = T_{D}(F)^{G} \subset
R_{D}(F)^{G} = E(F)[2]. $ \\
(3) \ The inverse images of $ 2 R_{D}(F), 2 E(F) $ under $
\varphi _{2},  \varphi _{1} $ respectively are given by \\
$ \varphi _{2}^{-1} (2 R_{D}(F)) = E(F) + R_{D}(F) = \varphi
_{1}^{-1} (2 E(F)). $
\par  \vskip 0.2cm

{\bf Lemma 1.3.} \\
$ (E(K) : E(F) + R_{D}(F)) = (T_{D}(F) : 2 R_{D}(F)) = (N(E(K)) : 2
E(F)). $
\par  \vskip 0.2cm

The proofs of Lemmas 1.1$\sim$1.3 are straightforward.
\par  \vskip 0.2cm

For the $ G-$modules $ E(K), E(F), R_{D}(F),  T_{D}(F) $ and $
N(E(K)), $ we have the following results about their corresponding
cohomology groups:
\par  \vskip 0.2cm

{\bf Proposition 1.4.} \ $ \text{H}^{1} (G, E (K)) = R_{D}(F) /
T_{D} (F); \quad \text{H}^{1} (G, E (F)) =  E (F)[2]; \\
\text{H}^{1}(G, R_{D}(F)) = R_{D}(F) / 2 R_{D}(F); \quad
\text{H}^{1}(G, T_{D}(F)) = T_{D}(F) / 2 T_{D}(F); \\
\text{H}^{1}(G, N(E(K))) = N(E(K))[2]. $
\par  \vskip 0.1cm
{\bf Proof.} \ Since $ G $ is cyclic, by the explicit formulae of
cohomology of finite cyclic groups (See [Se], pp.133, 128 for the
details), we have \\
$ \text{H}^{1} (G, E (K)) = \text{ker} \varphi _{1}
/ \text{im} \varphi _{2} = R_{D}(F) / T_{D} (F); \\
\text{H}^{1} (G, E (F)) = \text{ker}( \varphi _{1}|E(F)) /
\text{im}( \varphi _{2}|E(F)) = E (F)[2]; \\
\text{H}^{1}(G, R_{D}(F)) = \text{ker}( \varphi _{1}|R_{D}(F)) /
\text{im}( \varphi _{2}|R_{D}(F)) = R_{D}(F) / 2 R_{D}(F); \\
\text{H}^{1} (G, T_{D}(F)) = \text{ker}( \varphi _{1}|T_{D}(F)) /
\text{im}( \varphi _{2}|T_{D}(F)) = T_{D}(F) / 2 T_{D}(F); \\
\text{H}^{1} (G, N(E(K))) = \text{ker}( \varphi _{1}|N(E(K))) /
\text{im}( \varphi _{2}|N(E(K))) = N(E(K))[2]. $
\quad $ \Box $
\par  \vskip 0.2cm

{\bf Theorem 1.5.} \ The order of the group $ \text{H}^{1} (G, E
(K)) $ is
\begin{align*} \sharp \text{H}^{1} (G, E (K)) &= \frac{2^{r_{D, F}}
\cdot \sharp E(F)[2]}{(E(K) : E(F) + R_{D}(F))}
= \frac{2^{r_{D, F}} \cdot \sharp E(F)[2]}{(N(E(K)) : 2 E(F))} \\
&= 2^{r_{D, F} - r_{F}} \cdot (E(F) : N(E(K))).
\end{align*}
\par  \vskip 0.1cm
{\bf Proof.} \ Let $ A = R_{D}(F), B = E(K) $ and $ C = N(E(K)), $
their corresponding Herbrand quotients are $$ h(A) = h_{0}(A) /
h_{1}(A), \ h(B) = h_{0}(B) / h_{1}(B), \ h(C) = h_{0}(C) /
h_{1}(C), $$ where $ h_{m}( \cdot ) $ is the order of $
\widehat{\text{H}}^{m}(G, \cdot) \ (m = 0, 1) $ (see [AW, p.109]).
Since $ 2 E (F) \subset N(E(K)) \subset E(F), \ \text{rank}E_{D}(F)
= \text{rank} R_{D}(F) $ and $ \text{rank} N(E(K)) = \text{rank}
E(F), $ by Lemma 1.2 and Proposition 1.4, we have
\begin{align*}
h(R_{D}(F)) &= \frac{\sharp (R_{D}(F)^{G} / \varphi
_{1}(R_{D}(F)))}{\sharp \text{H}^{1} (G, R_{D}(F))} = \frac{\sharp
E(F)[2]}{\sharp (R_{D}(F) / 2 R_{D}(F))} = 2^{- r_{D, F}}, \\
h(E(K)) &= \frac{\sharp (E(K)^{G} / \varphi _{1}(E(K)))}{\sharp
\text{H}^{1} (G, E(K))} = \frac{\sharp ( E(F) / N(E(K)))}{\sharp
\text{H}^{1} (G, E(K))} \\
&= \frac{2^{r_{F}} \cdot \sharp E(F)[2]}{\sharp \text{H}^{1}(G,
E(K)) \cdot (N(E(K)) : 2 E(F))}, \\
h(N(E(K))) &= \frac{\sharp (N(E(K))^{G} / \varphi
_{1}(N(E(K))))}{\sharp \text{H}^{1} (G, N(E(K)))} = \frac{\sharp
(N(E(K)) / 2 N(E(K)))}{\sharp N(E(K))[2]} = 2^{r_{F}}.
\end{align*}
Since $ O \longrightarrow R_{D}(F) \longrightarrow E (K)
\longrightarrow ^{\varphi _{1}} N(E(K)) \longrightarrow O $ is an
exact sequence of $ G-$modules, by the theorem of Herbrand quotient
(see [AW], Prop.10 on p.109), we have $ h(E(K)) = h(R_{D}(F)) \cdot
h(N(E(K))). $ Therefore by the above calculation and Lemma 1.3, we
get
\begin{align*} \sharp \text{H}^{1} (G, E (K)) &= \frac{2^{r_{D, F}}
\cdot \sharp E(F)[2]}{(N(E(K)) : 2 E(F))} = \frac{2^{r_{D, F}}
\cdot \sharp E(F)[2]}{(E(K) : E(F) + R_{D}(F))} \\
&= 2^{r_{D, F} - r_{F}} \cdot (E(F) : N(E(K))).  \quad \quad \quad \Box
\end{align*}
\par  \vskip 0.2cm

{\bf Corollary 1.6.} \ If $ r_{F} = 0 $ and $ E(F)[2] = \{ O \}, $
then $ E(K) = E(F) + R_{D}(F) $ and $ \sharp \text{H}^{1} (G, E (K))
= 2^{r_{D, F}} = 2^{r_{K}}. $
\par  \vskip 0.1cm
{\bf Proof.} \ If $ r_{F} = 0 $ and $ E(F)[2] = \{ O \}, $ then by
the Mordell-Weil theorem, $ E(F) / 2 E(F) \cong (\Z / 2 \Z)^{r_{F}}
\oplus E(F)[2] = 0. $ So $ (N(E(K)) : 2 E(F)) = 1 $ because $
N(E(K)) / 2 E(F) \subset E(F) / 2 E(F), $ and then the conclusions
follow from Lemma 1.3 and Theorem 1.5.
\quad $ \Box $

\par     \vskip  0.4 cm

\hspace{-0.8cm}{\bf 2. A refined version of Yu's formula in the case
of elliptic curves}

\par \vskip 0.2 cm

For the quadratic extension $ K / F $ of number fields and the
elliptic curve $ E $ (over $ F $) as above, write $ M_{F} $ (resp. $
M_{K} $) for a complete set of places on $ F $ (resp. $ K $), let $
S_{\infty } $ be the set of infinite (i.e., Archimedean) places of $
F $ and $ S $ be the set of finite places of $ F $ obtained by
collecting together all places that ramify in $ K / F $ and all
places of bad reduction for $ E / F. $ Fix a place $ w \in M_{K} $
lying above $ v $ for each $ v \in M_{F}. $ Denote $
\text{Gal}(K_{w} / F_{v}) $ by $ G_{w}, $ where $ F_{v} $ and $
K_{w} $ are the completions of $ F $ at $ v $ and $ K $ at $ w, $
respectively. For each real place $ v \in S_{\infty }, $ let $
\sigma _{v}: F \rightarrow F_{v} = \R $ be the corresponding real
embedding, so $ \sigma _{v}(a) \in \R $ for any $ a \in F. $ For
each finite place $ v $ of $ F, $ we use $ v (\cdot ) $ to denote
the normalized additive valuation of $ F_{v}, $ i.e., $
v(F_{v}^{\ast}) = \Z. $ Let $
\parallel a \parallel _{F_{v}} = (\sharp k_{v})^{- v (a)} \ (a
\in F_{v}) $ denote the absolute value on $ F_{v} $ ($ k_{v} $ is
the residue field of $ F_{v} $), so is the meaning of $ \parallel a
\parallel _{K_{w}} $ on $ K_{w}. $ Let $ \Delta _{v}, \Delta _{D,
v}, $ and $ \Delta _{w} $ be the minimal discriminants for $ E $
over $ F_{v}, \ E_{D} $ over $ F_{v} $ and $ E $ over $ K_{w} $ (see
[Si1, p.186]), let $ c_{v}, c_{D, v} $ and $ c_{w} $ be the Tamagawa
numbers for $ E $ over $ F_{v}, \ E_{D} $ over $ F_{v} $ and $ E $
over $ K_{w} $ (see [Si1, p.451]), and let $ d (K_{w} / F_{v}) $ be
the discriminant of $ K_{w} / F_{v}, $ determined up to the square
of a unit of $ F_{v} $ (see [KT, p.332]). We also let $ ( \ , \
)_{F_{v}} $ denote the Hilbert norm-residue symbol, a
bimultiplicative form \ $ ( \ , \ )_{F_{v}}: \ F_{v}^{\ast } \times
F_{v}^{\ast }\rightarrow \mu _{2} = \{1, -1 \} $ whose properties
are described in [Se, Chapt.XIV]. For a vector space $ V $ over $
\F_{2}, $ the finite field with $ 2-$elements, we denote its
dimension by $ \dim _{2} V. $ Moreover, for
$ v \in M_{F} \setminus S_{\infty }, $
we denote $ \delta _{v} = \log _{2}(E(F_{v}) : N(E(K_{w}))). $ Then
the Theorem 7.6 in [KT, p.332] states (in an equivalent form) that
$$ \delta _{v} = \log _{2} \left( \frac{c_{v} c_{D, v}}{c_{w}} \left(
\frac{\parallel \Delta _{v} \Delta _{D, v} d (K_{w} / F_{v})^{-6}
\parallel _{F_{v}}}{\parallel \Delta _{w} \parallel
_{K_{w}}}\right)^{1 / 12} \right). $$

{\bf Theorem 2.1 (A refined version of Yu's formula).} \ Assume
that the Shafarevich-Tate groups $ \amalg\hskip-6pt\amalg(E/F),
\amalg\hskip-6pt\amalg(E_{D}/F) $ and $ \amalg\hskip-6pt\amalg(E/K) $
are finite. Then
$$ \frac{\sharp \amalg\hskip-6pt\amalg(E/F) \cdot \sharp
\amalg\hskip-6pt\amalg(E_{D}/F)}{\sharp \amalg\hskip-6pt\amalg(E/K)}
= 2^{r_{D, F} - r_{F} - \delta (E, F, K)} \cdot (E(F) :
N(E(K)))^{2}, $$ where \ $ \delta (E, F, K) = \delta _{\infty } +
\delta _{f}, $ \ with $ \delta _{\infty } = \sharp \{v \in S_{\infty
} : \ v \ \text{is ramified in} \ K \ \text{and} \ \sigma _{v} (
\Delta (E)) > 0 \} $ ($ \Delta (E) $ is the discriminant of the
elliptic curve $ E $ over $ F $), \ and $ \delta _{f} = \sum _{v
\in S_{0}} \delta _{v}. $ Here $ S_{0} = \{v \in S: \ v
\ \text{is ramified or inertial in} \ K \}. $
\par  \vskip 0.2cm

{\bf Proof.} \ By Yu's formula [Y] $$ \frac{\sharp
\amalg\hskip-6pt\amalg(E/F) \cdot \sharp
\amalg\hskip-6pt\amalg(E_{D}/F)}{\sharp \amalg\hskip-6pt\amalg(E/K)}
= \frac{\sharp \widehat{\text{H}}^{0}(G, E(K)) \cdot \sharp
\text{H}^{1}(G, E(K))}{\prod _{v \in M_{F}} \sharp \text{H}^{1}
(G_{w}, E(K_{w}))}. $$ By definition, $ \widehat{\text{H}}^{0}(G,
E(K)) = E(K)^{G} / (1 + \sigma )E(K) = E(F) / N(E(K)), $ so by the
above Theorem 1.5, we get $$ \frac{\sharp
\amalg\hskip-6pt\amalg(E/F) \cdot \sharp
\amalg\hskip-6pt\amalg(E_{D}/F)}{\sharp \amalg\hskip-6pt\amalg(E/K)}
= \frac{2^{r_{D, F} - r_{F}} \cdot (E(F) : N(E(K)))^{2}}{\prod _{v
\in M_{F}} \sharp \text{H}^{1} (G_{w}, E(K_{w}))}. $$ On the other
hand, by the Corollary 4.4 in [Ma, p.204], we have $ \text{H}^{1}
(G_{w}, E(K_{w})) = 0 $ for any $ v \notin S \cup S_{\infty }. $
Therefore
$$ \frac{\sharp \amalg\hskip-6pt\amalg(E/F) \cdot \sharp
\amalg\hskip-6pt\amalg(E_{D}/F)}{\sharp \amalg\hskip-6pt\amalg(E/K)}
= \frac{2^{r_{D, F} - r_{F}} \cdot (E(F) : N(E(K)))^{2}}{\prod _{v
\in S \cup S_{\infty }} \sharp \text{H}^{1} (G_{w}, E(K_{w}))}.
\eqno(2.1) $$ By our assumption, the Shafarevich-Tate groups are
finite; also $ (E(F) : N(E(K))) < \infty $ because $
\text{rank}E(F) = \text{rank}N(E(K)), $ so by the above formula
(2.1), $ \text{H}^{1} (G_{w}, E(K_{w})) $ is a finite set
for each $ v \in S \cup S_{\infty }. $ \\
Let $ v \in S_{\infty }, $ if $ v $ is unramified in $ K, $ then $
\text{H}^{1} (G_{w}, E(K_{w})) = 0 $ because $ K_{w} = F_{v} = \R $
or $ \C. $ So we may assume that $ v $
is ramified in $ K, $ then $ F_{v} = \R $ and $ K_{w} = \C. $ By the
Theorem 2.4 of Chapter V in [Si2], we have
$$ \text{H}^{1} (G_{w}, E(K_{w})) = \text{H}^{1}
(\text{Gal}(\C / \R), E(\C)) \cong \left \{
\begin{array}{l}
0 \ \quad \ \text{if} \ \sigma _{v} ( \Delta (E)) < 0, \\ \\
\Z / 2 \Z \quad \text{if} \ \sigma _{v} ( \Delta (E)) > 0.
\end{array}
\right. $$ Hence $$ \prod _{v \in S_{\infty }} \sharp \text{H}^{1}
(G_{w}, E(K_{w})) = \sharp \left( \Z / 2 \Z \right)^{\delta _{\infty
}} = 2^{\delta _{\infty }}. \eqno(2.2) $$ Let $ v \in S, $ if $ v
\notin S_{0}, $ then $ v $ splits completely in $ K, $ so $ K_{w} =
F_{v} $ and then $ \text{H}^{1} (G_{w}, E(K_{w})) = 0. $ For $ v \in
S_{0}, $ since $ \text{H}^{1} (G_{w}, E(K_{w})) $ is finite as
mentioned above, by Tate local duality [Ta1] (see also [Ma,
Prop.4.2]), we have $ \sharp \text{H}^{1} (G_{w}, E(K_{w})) =
(E(F_{v}) : N(E(K_{w}))). $ Hence by the Theorem 7.6 and the Remark in
[KT, pp. 332, 333] (or by Prop.$ 1\sim 5 $ in [Kr]), we get
$$ \prod _{v \in S} \sharp \text{H}^{1} (G_{w}, E(K_{w}))
= \prod _{v \in S_{0}} \sharp \text{H}^{1} (G_{w}, E(K_{w})) =
2^{\delta _{f}}. \eqno(2.3)
$$ Substitute (2.2) and (2.3) into (2.1), we get
$$ \frac{\sharp \amalg\hskip-6pt\amalg(E/F) \cdot \sharp
\amalg\hskip-6pt\amalg(E_{D}/F)}{\sharp \amalg\hskip-6pt\amalg(E/K)}
= 2^{r_{D, F} - r_{F} - \delta (E, F, K)} \cdot (E(F) :
N(E(K)))^{2}. $$ The proof is completed. \quad $ \Box $
\par  \vskip 0.2cm

{\bf Remark.} \ By the results of Kramer on the local norm index in
[Kr], one can calculate $ \delta _{f} $ (hence $ \delta (E, F, K) $)
for most cases as follows: \\
$ \delta _{f} = \delta _{g} + \delta _{m} + \delta _{a}, $ where $
\delta _{g}, \delta _{m} $ and $ \delta _{a} $ are defined as
follows:
\begin{align*}
&\delta _{a} = \sum _{v \in S_{a}} \delta _{v}; \\
&\delta _{m} = \delta _{smr} + \delta _{nsmr} \ \text{with} \ \delta
_{smr} = \frac{1}{2} \sum _{v \in S_{smr}} \left(1 + (\Delta _{v},
D)_{F_{v}} \right) \ \text{and} \\
&\delta _{nsmr} = \frac{1}{2} \sum _{v \in S_{nsmr}^{\prime }}
\left(1 + (-1)^{v (\Delta _{v})} \right) + \sum _{v \in
S_{nsmr}^{\prime \prime}} \left( \frac{1}{2} \left(1 + (\Delta _{v},
D)_{F_{v}} \right) \cdot (-1)^{v (\Delta _{v})} + 1 \right); \\
&\delta _{g} = \sum _{v \in S_{g}} \dim _{2}
\widetilde{E_{v}}(k_{v}) [2] + \sum _{v \in S_{gu}} \varepsilon (v), \quad 
\text{where}
\end{align*}
 $$ \varepsilon (v) = \left \{\begin{array}{l}
\frac{1}{2} \left(1 - (-1)^{v(D)} \right) \cdot [F_{v} : \Q_{2}] \
\text{if} \ E \ \text{has good supersingular reduction at} \ v, \\
\\
\frac{1}{2} (3 + (\Delta _{v}, D)_{F_{v}}) \quad \quad \text{if} \ E
\ \text{has good ordinary reduction at} \ v.
\end{array}
\right. $$ Here $ \widetilde{E_{v}} $ is the reduction of $ E $ at $
v, \ k_{v} $ is the residue field of $ F_{v}, $ \\
$ S_{0} = \{v \in S: \ v \ \text{is ramified or inertial in} \ K \}; \\
S_{g} = \{v \in S_{0}: \ v \nmid 2 \ \text{and} \ E \
\text{has good reduction at} \ v \}; \\
S_{gu} = \{v \in S_{0}: \ v \mid 2, \ E \ \text{has good reduction
at} \ v \ \text{and} \ F_{v} \ \text{is unramified over} \
\Q_{2} \}; \\
S_{ar} = \{v \in S_{0}: \ E \ \text{has additive reduction at} \ v
\}; \\
S_{a} = S_{ar} \cup \{v \in S_{0}: \ v \mid 2, \ E \ \text{has good
reduction at} \ v \ \text{and} \ F_{v} \
\text{is ramified over} \ \Q_{2} \}; \\
S_{smr} = \{v \in S_{0}: \ E \ \text{has split multiplicative
reduction at} \ v \}; \\
S_{nsmr} = \{v \in S_{0}: \ E \ \text{has non-split
multiplicative reduction at} \ v \} \\
= S_{nsmr}^{\prime } \sqcup S_{nsmr}^{\prime \prime} \ (\text{the
disjoint union}), \quad
\text{where} \\
S_{nsmr}^{\prime } = \{v \in S_{nsmr}: \ v \ \text{is inertial
in} \ K \}, \\
S_{nsmr}^{\prime \prime} = \{v \in S_{nsmr}: \ v \ \text{is
ramified in} \ K \}. $ \\
Obviously, $ S_{0} = S_{g} \sqcup S_{gu} \sqcup S_{a} \sqcup S_{smr}
\sqcup S_{nsmr} $ (the disjoint union).
\par     \vskip  0.2 cm

{\bf A Note Added.} \ The main aim of this paper is to explicitly
work out the orders of Shafarevich-Tate groups of some elliptic curves
in quadratic number fields by refining Yu's formula. Nevertheless, as
pointed out by the anonymous referee, the above methods can apply to more
general cases. More precisely, as stated in the referee's comments, the
above Theorem 1.5 holds verbatim for an abelian variety over a
number field $ F $ (or even a global field of characteristic
different from $ 2 $), and, except for Lemma 1.1, the entire
discussion can take place in the context of a group $ G $ of order $
2 $ acting linearly on a finitely generated abelian group $ A. $
Then the refined Yu's formula for an abelian variety $ A / F $ is
$$ \frac{\sharp \amalg\hskip-6pt\amalg(A/F) \cdot \sharp
\amalg\hskip-6pt\amalg(A_{D}/F)}{\sharp \amalg\hskip-6pt\amalg(A/K)}
= \frac{2^{r_{D, F} - r_{F}} \cdot (A^{\prime }(F) :
N(A^{\prime}(K))) \cdot (A(F) : N(A(K)))}{\prod _{v}(A(F_{v}) : N(A(K_{w})))},
$$ where $ A^{\prime } $ is the dual abelian variety of $ A. $ \\
It would be interesting to find a similar formula of Kramer-Tunnell 
for higher-dimensional abelian varieties to make these local norm 
cokernel quantities (i.e. $(A(F_{v}) : N(A(K_{w})))$) explicitly calculated.

\par     \vskip  0.4 cm

\hspace{-0.8cm}{\bf 3. Application I - congruent numbers elliptic
curves}
\par \vskip 0.2 cm

Let $ n \in \Z \setminus \{0, 1 \} $ be a square free integer and $
K = \Q (\sqrt{n}) $ be a quadratic number field. In this section, we
consider elliptic curves $ E : y^{2} = x^{3} - x $ and $ E_{n} :
y^{2} = x^{3} - n^{2} x. $ All these curves have complex
multiplication by $ \Z[\sqrt{-1}], $ the Gaussian integral ring. Let
$ w \in M_{K} $ be a place of $ K $ lying over $ 2, $ as in section
2 above, recall that the notations $ \Delta _{w} $ and $ c_{w} $
represent the minimal discriminant and the Tamagawa factor for $ E $
over $ K_{w}, $ respectively. Denote by $ \text{ord}_{w} $ the
normalized additive valuation of $ K_{w}. $
\par \vskip 0.2 cm

{\bf Lemma 3.1.} \ We have
$$ \text{ord}_{w}(\Delta _{w}) = \left \{\begin{array}{l} 6
\quad \quad \ \text{if} \ n \equiv 1 \ (\text{mod} \ 4), \\
12 \quad \quad  \text{if} \ n \equiv 2 \ \text{or} \ 3 \ (\text{mod}
\ 4), \quad \text{and} \end{array} \right. $$
$$ c_{w} = \left \{\begin{array}{l} 4 \quad
\quad \text{if} \ n \equiv 2 \ \text{or} \ 7 \ (\text{mod} \ 8), \\
2 \quad \quad \text{if} \ n \equiv 1, 3, 5 \ \text{or} \ 6 \
(\text{mod} \ 8).
\end{array} \right. $$
\par  \vskip 0.1cm
{\bf Proof.} \ Up to isomorphisms, there are
exactly seven quadratic extensions of $ \Q_{2}, $ namely, $
\Q_{2}(\sqrt{-1}), \Q_{2}(\sqrt{-2}), \Q_{2}(\sqrt{2}),
\Q_{2}(\sqrt{-3}), \Q_{2}(\sqrt{3}), 
\Q_{2}(\sqrt{-6}), \Q_{2}(\sqrt{6}) $ (see [W, p.248]).
Furthermore, one can easily verify that \\
$ K_{w} \cong \Q_{2} \Longleftrightarrow n \equiv 1 \ (\text{mod} \
8); \quad K_{w} \cong \Q_{2}(\sqrt{-3}) \Longleftrightarrow n \equiv
5 \ (\text{mod} \ 8); $ \\
$ K_{w} \cong \Q_{2}(\sqrt{-1}) \Longleftrightarrow n \equiv 7 \
(\text{mod} \ 8); \quad K_{w} \cong \Q_{2}(\sqrt{3})
\Longleftrightarrow n \equiv 3 \ (\text{mod} \ 8); $ \\
$ K_{w} \cong \Q_{2}(\sqrt{-2}) \Longleftrightarrow n \equiv 14 \
(\text{mod} \ 16); \quad K_{w} \cong \Q_{2}(\sqrt{2})
\Longleftrightarrow n \equiv 2 \ (\text{mod} \ 16); $ \\
$ K_{w} \cong \Q_{2}(\sqrt{-6}) \Longleftrightarrow n \equiv 10, 26
\ \text{or} \ 42 \ (\text{mod} \ 48); $ \\
$  K_{w} \cong \Q_{2}(\sqrt{6}) \Longleftrightarrow n \equiv 6, 22 \
\text{or} \ 38 \ (\text{mod} \ 48). $ \\
Next, by Tate's algorithm (see [Ta2, pp.47$\thicksim$52], [Si2,
Chapt.IV$\S9$]), we get $$ c_{w} = \left \{\begin{array}{l} 2 \quad
\quad \text{if} \ K_{w} = \Q_{2}, \Q_{2}(\sqrt{-3}),
\Q_{2}(\sqrt{-2}),
\Q_{2}(\sqrt{3}) \ \text{or} \ \Q_{2}(\sqrt{6}), \\
4 \quad \quad \text{if} \ K_{w} = \Q_{2}(\sqrt{-1}),
\Q_{2}(\sqrt{2}) \ \text{or} \ \Q_{2}(\sqrt{-6}), \quad \text{and}
\end{array} \right. $$
$$ \text{ord}_{w}(\Delta _{w}) = \left \{\begin{array}{l} 6
\quad \quad \ \text{if} \ n \equiv 1 \ (\text{mod} \ 4), \\
12 \quad \quad  \text{if} \ n \equiv 2 \ \text{or} \ 3 \ (\text{mod}
\ 4), \end{array} \right. $$ from which the conclusion follows, and
the proof is completed. \quad $ \Box $
\par \vskip 0.2 cm

{\bf Lemma 3.2.} \ We have
$$ \delta (E, \Q, K) =
\left \{\begin{array}{l} 2 \omega _{0}(n) \quad \quad \quad \ \
\text{if} \ n > 0 \ \text{and} \ n \equiv 1 \ (\text{mod} \ 8), \\
1 + 2 \omega _{0}(n) \quad \quad \text{if} \ n > 0 \ \text{and} \
n \equiv 5 \ \text{or} \ 7 \ (\text{mod} \ 8), \\
3 + 2 \omega _{0}(n) \quad \quad \text{if} \ n > 0 \ \text{and} \
n \equiv 6 \ (\text{mod} \ 8), \\
2 + 2 \omega _{0}(n) \quad \quad \text{if} \ n > 0 \ \text{and} \
n \equiv 2 \ \text{or} \ 3 \ (\text{mod} \ 8), \\
1 + 2 \omega _{0}(n) \quad \quad  \text{if} \ n < 0 \ \text{and} \
n \equiv 1 \ (\text{mod} \ 8), \\
2 + 2 \omega _{0}(n) \quad \quad \text{if} \ n < 0 \ \text{and} \
n \equiv 5 \ \text{or} \ 7 \ (\text{mod} \ 8), \\
3 + 2 \omega _{0}(n) \quad \quad \text{if} \ n < 0 \ \text{and} \
n \equiv 2 \ \text{or} \ 3 \ (\text{mod} \ 8), \\
4 + 2 \omega _{0}(n) \quad \quad \text{if} \ n < 0 \ \text{and} \ n
\equiv 6 \ (\text{mod} \ 8),
\end{array} \right. $$
where $ \omega _{0}(n) $ is the number of odd prime divisors of $ n.
$
\par  \vskip 0.1cm
{\bf Proof.} \ Since $ \Delta (E) = 64 > 0, \ E $ has good reduction
everywhere except at $ 2 $ with additive reduction. So, by
definition, $ \ S = \{2 \} \cup \{p : \ p \ \text{is a prime and} \
p \mid n \}, S_{gu} = S_{smr} = S_{nsmr} = \emptyset $ and $ S_{g} =
S \setminus \{ 2 \}. $ So $ \delta _{m} = 0, $ and $ \delta _{\infty
} = 0 $ (resp., $ 1 $) if $ n > 0 $ (resp., $ n < 0 $). Moreover,
for each odd prime $ p, \ E $ has good reduction at $ p, $ and it
easy to see that $ \widetilde{E}(\F_{p})[2] \cong \left(\Z / 2 \Z
\right)^{2}, $ so $ \delta _{g} = \sum _{p \in S_{g}} \dim _{2}
\widetilde{E}(\F_{p})[2] = 2 \omega _{0}(n). $ Hence by definition,
$ \delta (E, \Q, K) = \delta _{\infty } + \delta _{g} + \delta _{m}
+ \delta _{a} = 2 \omega _{0}(n) + \delta _{\infty } + \delta _{a}.
$ We divide our
discussion into the following cases. \\
Case A. $ n \equiv 1 \ (\text{mod} \ 8). $ Then $ 2 $ splits
completely in $ K, $ and then $ S_{a} = \emptyset , $ so $ \delta
_{a} = 0, $ which implies $ \delta (E, \Q, K) = 2 \omega _{0}(n) $
(resp., $ 2 \omega _{0}(n) + 1 $ ) if $ n > 0 $ (resp., $ n < 0 $). \\
Case B. $ n \equiv 2, 3, 5, 6 \ \text{or} \ 7 \ (\text{mod} \ 8). $
Then $ 2 $ is ramified or inertial in $ K, $ so $ S_{a} = \{ 2 \}. $
Let $ w \in M_{K} $ be the unique place in $ K $ lying above $ 2, $
then $ K_{w} = \Q_{2}(\sqrt{n}) $ is a quadratic extension over $
\Q_{2}. $ By Thm.7.6 in [KT], we get
$$ \delta _{a} = \delta _{2} = \log _{2}(E(\Q_{2}) : N(E(K_{w})))
= \log _{2} \left( \frac{c_{2} c_{n, 2}}{c_{w}} 
\left( \frac{\parallel \Delta _{2} \Delta _{n, v} d_{w}^{-6}
\parallel _{\Q_{2}}}{\parallel \Delta _{w} \parallel
_{K_{w}}} \right)^{1 / 12} \right). \eqno(3.1) $$ Now we only need
to compute all the values of $ c_{2}, c_{n, 2}, c_{w}, \Delta _{2},
\Delta _{n, v}, \Delta _{w} $ and $ d_{w}. $ Firstly, by a method in
([KT], p.331)
$$ d_{w} = d (K_{w} / \Q_{2}) =
\left \{\begin{array}{l} n \quad \quad \
\text{if} \ n \equiv 5 \ (\text{mod} \ 8), \\
4 n \quad \quad \text{if} \ n \equiv 2 \ \text{or} \ 3 \ (\text{mod}
\ 4).
\end{array} \right. $$
Next, for the elliptic curves $ E $ and $ E_{n} $ over $ \Q_{2}, $
by Tate's algorithm (see [Ta2, pp.47$\thicksim$52], [Si2,
Chapt.IV$\S9$]), one can easily obtain
that $ v_{2}(\Delta _{2}) = 6, \ c_{2} = 2 $ and \\
$ v_{2}(\Delta _{n, 2}) = 6 $ if $ n \equiv 3, 5 \ \text{or} \ 7 \
(\text{mod} \ 8);  \quad  v_{2}(\Delta _{n, 2}) = 12 $ if $ n \equiv
2 \ \text{or} \ 6 \ (\text{mod} \ 8); $ \\
$ c_{n, 2} = 2 $ if $ n \equiv 3, 5 \ \text{or} \ 7 \ (\text{mod} \
8); \quad c_{n, 2} = 4 $ if $ n \equiv 2 \ \text{or} \ 6 \
(\text{mod} \ 8). $ \\
Also by Lemma 3.1 above, we have \\
$ \text{ord}_{w}(\Delta _{w}) = 12 $ if $ n \equiv 2 \ \text{or} \ 3
\ (\text{mod} \ 4);  \quad \text{ord}_{w}(\Delta _{w}) = 6 $ if $ n
\equiv 5 \ (\text{mod} \ 8); $ \\
$ c_{w} = 2 $ if $ n \equiv 3, 5 \ \text{or} \ 6 \ (\text{mod} \ 8);
\quad c_{w} = 4 $ if $ n \equiv 2 \ \text{or} \ 7 \
(\text{mod} \ 8). $ \\
Now substitute all of them into (3.1), the conclusion for case B
then follows, and the proof is completed. \quad $ \Box $
\par \vskip 0.2 cm

It is well known that the $ L-$ function $ L(E/\Q, s) = \sum b_{m}
m^{-s} $ of the elliptic curve $ E = E_{1} : y^{2} = x^{3} - x $
corresponds to a weight two cusp form $ g = \sum b_{m} q^{m} \in
S_{2}(\Gamma _{0} (32)) $ (see [Kob, p.217]), and for the elliptic
curve $ E_{n} : y^{2} = x^{3} - n^{2} x, $ by Tunnell's theorem (see
[T, p.328] or [Kob, p.217]), there exist a form $ f = \sum a_{m}
q^{m} \in S_{3/2}(\widetilde{\Gamma} _{0} (128)) $ and a form $
f^{\prime } = \sum a_{m}^{\prime } q^{m} \in
S_{3/2}(\widetilde{\Gamma} _{0} (128), \chi _{2}) $ such that their
Shimura lifts $ \text{Shimura} (f) = \text{Shimura} (f^{\prime }) =
g $ and $$ L(E_{n}/\Q, 1) = \left \{\begin{array}{l} \frac{\omega
}{4 \sqrt{n}} a_{n}^{2} \quad
\text{if} \ n \ \text{is odd,} \\ \\
\frac{\omega }{2 \sqrt{n}} (a^{\prime }_{n /2})^{2} \quad \text{if}
\ n \ \text{is even.}
\end{array} \right. $$
where $ \omega = \int _{1}^{\infty} \frac{dx}{\sqrt{x^{3} - x}} =
2.6220575 $ is the least positive period of $ E / \Q. $
\par \vskip 0.2 cm

{\bf Theorem 3.3.} \ Let $ n $ be a square
free integer satisfying one of the following conditions \\
(1) \ $ n > 0 $ and $ n \equiv 1, 2 \ \text{or} \ 3 \ (\text{mod} \
8); $ \quad (2) \ $ n < 0 $ and $ n \equiv 5, 6 \ \text{or} \ 7 \
(\text{mod} \ 8). $ \\
Then for the elliptic curves $ E_{n}: y^{2} = x^{3} - n^{2} x $ and
$ E = E_{1} $ as above, if the full BSD conjecture (see [Si1,
p.452]) is true for $ E_{n} $ over $ \Q $ with $ L(E_{n}/\Q, 1) \neq
0, $ and $ \amalg\hskip-6pt\amalg(E/\Q(\sqrt{n})) $ is finite, we
have
$$ \sharp \amalg\hskip-6pt\amalg(E/\Q(\sqrt{n})) =
\left \{\begin{array}{l} 2^{-4} \cdot a_{n}^{2} \quad \quad \
\text{if} \ n > 0 \ \text{and} \ n
\equiv 1 \ (\text{mod} \ 8), \\ \\
2^{-2} \cdot a_{n}^{2} \quad \quad \ \text{if} \ n > 0 \ \text{and}
\ n \equiv 3 \ (\text{mod} \ 8), \\ \\
2^{-2} \cdot (a^{\prime }_{n/2})^{2} \quad \text{if} \ n > 0 \
\text{and} \ n \equiv 2 \ (\text{mod} \ 8), \\ \\
2^{-2} \cdot a_{- n}^{2} \quad \quad \text{if} \ n < 0 \ \text{and}
\ n \equiv 5 \ \text{or} \ 7 \ (\text{mod} \ 8), \\  \\
(a^{\prime }_{- n/2})^{2} \quad \quad \ \text{if} \ n < 0 \
\text{and} \ n \equiv 6 \ (\text{mod} \ 8),
\end{array} \right. $$
where $ a_{\mid n \mid } $ and $ a^{\prime }_{\mid n/2 \mid } $ are
the Fourier coefficients of the above modular forms $ f $ and $
f^{\prime }. $
\par  \vskip 0.1cm
{\bf Proof.} \ We prove the case that $ n > 0 $ satisfying $ n
\equiv 1 \ (\text{mod} \ 8), $ the other cases can be similarly
verified. For this case, by Lemma 3.2 above,
$ \delta (E, \Q, K) = 2 \omega _{0}(n). $ By the assumption, $
L(E_{n}/\Q, 1) \neq 0 $ and the full BSD conjecture is true for $
E_{n} $ over $ \Q, $ so $ r_{n, \Q} = 0, \
\amalg\hskip-6pt\amalg(E_{n}/\Q) $ is finite and $ L(E_{n}/\Q, 1) /
\Omega _{E_{n}/\Q} = \text{(BSD)}_{\infty , \Q} (E_{n}) $ with $
\Omega _{E_{n}/\Q} = \omega / \sqrt{n}, $ where
$$ \text{(BSD)}_{\infty , \Q} (E_{n}) =
\text{Reg}_{\infty , \Q}(E_{n}) \times \frac{\sharp
\amalg\hskip-6pt\amalg(E_{n}/\Q) \prod _{v \in M_{\Q}} c_{v}
}{\sqrt{d(\Q)} \times \sharp E_{n}(\Q)_{\text{tors}}^{2}}. $$ We
have $ \text{Reg}_{\infty , \Q}(E_{n}) = 1 $ because $ r_{n, \Q} =
0; $ obviously, $ d(\Q) = 1; $ also $ \sharp E_{n}(\Q)_{\text{tors}}
= 4 $ (see [Si1], pp.346, 347); moreover, $ c_{\infty } = 2 $
because $ E_{n} $ is not connected over $ \R, $ then from [R2,
p.235] we have
\begin{align*} & \prod _{v \in M_{\Q}} c_{v} = c_{\infty }
\cdot \prod _{p < \infty } c_{p} = 2 \times 2^{2 \omega _{0}(n) + 1}
= 2^{2 \omega _{0}(n) + 2}, \quad \text{hence} \\
&L(E_{n}/\Q, 1) = \Omega _{E_{n}/\Q} \times \text{(BSD)}_{\infty ,
\Q} (E_{n}) = 2^{2 \omega _{0}(n) - 2} \times \frac{\omega
}{\sqrt{n}} \times \sharp \amalg\hskip-6pt\amalg(E_{n}/\Q). \quad
(3.2)
\end{align*}
On the other hand, by Tunnell's theorem (see [T, Thm.3, p.328] or
[Kob, p.217]), we have $ L(E_{n}/\Q, 1) = \omega a_{n}^{2} /
(4\sqrt{n}) $ with the Fourier coefficient $ a_{n} $ of the modular
form $ f $ mentioned above. Therefore by (3.2), we get $ \sharp
\amalg\hskip-6pt\amalg(E_{n}/\Q)
= 2^{-2 \omega _{0}(n)} \cdot a_{n}^{2}. $ As mentioned before, 
$ \amalg\hskip-6pt\amalg(E/\Q) = 0, $ and $ E(\Q) = E(\Q)[2] \cong
\left(\Z / 2 \Z \right)^{2}, $ so $ r_{K} = r_{n, \Q} + r_{\Q} = 0,
$ and it is easy to know that $ E(K)_{tors} = E(\Q)[2], $ hence by
definition, we have $ (E(\Q) : N(E(K))) = 4. $ By assumption, $
\amalg\hskip-6pt\amalg(E/\Q(\sqrt{n})) $ is finite, hence by the
above refined Yu's formula, we get
\begin{align*} \sharp \amalg\hskip-6pt\amalg(E/\Q(\sqrt{n})) &= 2^{-r_{n,
\Q} + r_{\Q} + \delta (E, \Q, K)} \cdot (E(\Q) : N(E(K)))^{-2}
\cdot \sharp \amalg\hskip-6pt\amalg(E/\Q) \cdot \sharp
\amalg\hskip-6pt\amalg(E_{n}/\Q) \\
&= 2^{2 \omega _{0}(n)} \cdot 4^{-2} \cdot 2^{-2 \omega _{0}(n)}
\cdot a_{n}^{2} = 2^{-4} a_{n}^{2}. \quad \text{In particular,} \\
&\frac{\sharp \amalg\hskip-6pt\amalg(E/\Q(\sqrt{n}))}{\sharp
\amalg\hskip-6pt\amalg(E_{n}/\Q)} = 2^{2 \omega _{0}(n) - 4}.
\end{align*}
Therefore the conclusion of the Shafarevich-Tate groups $
\amalg\hskip-6pt\amalg(E/\Q(\sqrt{n})) $ is obtained,
and the proof of is completed.
\quad $ \Box $
\par \vskip 0.2 cm

{\bf Remark.} \ (1) \ Note that $ E $ and $ E_{n} $ are isomorphic
over $ \Q(\sqrt{n}), $ and $ E_{n} = E_{- n}, $ so in particular $
\sharp \amalg\hskip-6pt\amalg(E_{n}/\Q(\sqrt{\pm n})) = \sharp
\amalg\hskip-6pt\amalg(E/\Q(\sqrt{\pm n})), $ and one has all the
same results for $ E_{n} $ over $ \Q(\sqrt{n}) $ as $ E $ in the
above Theorem 3.3. Moreover, the Fourier coefficients $ a_{n} $ and
$ a_{n}^{\prime } \ (n > 0) $ can be determined by the number of
solutions of some concrete quadratic forms in three variables (see
Tunnell's theorem in [T,p.323 and p.325] for the detail). \\
(2) \ The above proof of Theorem 3.3 also shows that the ratio \\
$ \sharp \amalg\hskip-6pt\amalg(E/\Q(\sqrt{n})) / \sharp
\amalg\hskip-6pt\amalg(E_{n}/\Q) $ can be arbitrarily large. For
example, for a square free positive integer $ n $ satisfying the condition
of the following Corollary 3.4.(3), then $ \sharp
\amalg\hskip-6pt\amalg(E/\Q(\sqrt{n})) / \sharp
\amalg\hskip-6pt\amalg(E_{n}/\Q) = (2^{-2} \cdot a_{n}^{2}) /
(2^{-2 \omega _{0}(n)} \cdot a_{n}^{2})= 2^{2 \omega _{0}(n) - 2}, $ where
$ \omega _{0}(n) $ is the number of odd prime divisors of $ n. $
\par \vskip 0.2 cm

{\bf Corollary 3.4.} \ For the elliptic curves $ E_{n} : y^{2} =
x^{3} - n^{2} x $ and $ E = E_{1} $ as above, assume that
$ \amalg\hskip-6pt\amalg(E/\Q(\sqrt{n})) $ is finite.  \\
(1) \ If $ n = \pm p, \ p $ is a prime number, and $ p \equiv 3 \
(\text{mod} \ 8), $ then $$ \sharp
\amalg\hskip-6pt\amalg(E/\Q(\sqrt{p})) = \sharp
\amalg\hskip-6pt\amalg(E/\Q(\sqrt{-p})) = \frac{1}{4} a_{p}^{2}, $$
in particular, all such $ a_{p} $ are even. \\
(2) \ Suppose $ n = \pm p_{1} \cdots p_{m} \equiv 1 \ (\text{mod} \
4), $ where $ p_{1}, \cdots , p_{m} $ are distinct prime numbers
with $ p_{i} \not\equiv 5 \ (\text{mod} \ 8). $ If $ s_{k} (n) = 1,
$ then
$$ \sharp \amalg\hskip-6pt\amalg(E/\Q(\sqrt{n})) =
\left \{\begin{array}{l} 2^{-4} \cdot a_{n}^{2} \quad \quad \
\text{if} \ n > 0 \ \text{and} \ n
\equiv 1 \ (\text{mod} \ 8), \\ \\
2^{-2} \cdot a_{-n}^{2} \quad \quad \ \text{if} \ n < 0 \ \text{and}
\ n \equiv 5 \ (\text{mod} \ 8).
\end{array} \right. $$
(3) \ Suppose $ n = p_{1} \cdots p_{m}, $ where $ p_{1}, \cdots ,
p_{m} $ are distinct prime numbers with $ p_{1} \equiv 3 \
(\text{mod} \ 8) $ and $ p_{2} \equiv \cdots \equiv p_{m} \equiv 1 \
(\text{mod} \ 8). $ If $ s_{2m -1} (- n) = 1, $ then $$ \sharp
\amalg\hskip-6pt\amalg(E/\Q(\sqrt{n})) = 2^{-2} \cdot a_{n}^{2}, \
\text{in particular,} \ 2^{m} \parallel a_{n}, \ i.e., \ v_{2}
(a_{n}) = m. $$ Here $ s_{k}(n) $ and $
s_{2m -1}(-n) $ are the $ \F_{2}-$valued functions on $ n $ and its
Gaussian prime factors defined in [Z, p.387].
\par \vskip 0.1 cm
{\bf Proof.} \ (1). \ By a theorem of Rubin (see [R1], P.26), the
full BSD conjecture is true for $ E_{p} : y^{2} = x^{3} - p^{2} x $
over $ \Q $ and $ L(E_{p}/ \Q, 1) \neq 0, $ so the conclusion
follows directly from the above Theorem 3.3. \\
(2) and (3). \ By the Theorem 2 and Proposition 3 in [Z, p.387], the
full BSD conjecture is true for $ E_{n} : y^{2} = x^{3} - n^{2} x $
over $ \Q $ and $ L(E_{n}/ \Q, 1) \neq 0 $ in these cases, so the
conclusion of the orders of the Shafarevich-Tate groups for $ E $ 
over $ \Q(\sqrt{n})$ follows directly from the above Theorem 3.3. 
Now we come to compute the $ 2-$adic
valuation of $ a_{n} $ in case (3). In fact, by the Theorem 2 in
[Z, p.387], we know that, if $ s_{2m -1} (- n) = 1, $ then $ L(E_{n}/\Q, 1)
\neq 0 $ and the $ 2-$Selmer group $ S^{(2)} (E_{n}/ \Q) $ has order
$ 4. $ So by a theorem of Coates-Wiles (see [CW, Thm.1, p.223]), one
has $ r_{n, \Q} = \text{rank} E_{n}(\Q) = 0, $ so $ E_{n}(\Q) / 2
E_{n}(\Q) \cong E_{n}(\Q)[2] \cong \left( \Z / 2 \Z \right)^{2}. $
Then by the exact sequence (see [Si1], chapt.X, Thm.4.2)
$$ 0 \rightarrow E_{n}(\Q) / 2 E_{n}(\Q) \rightarrow S^{(2)} (E_{n}/ \Q)
\rightarrow \amalg\hskip-6pt\amalg(E_{n}/\Q)[2] \rightarrow 0 $$ we
get $ \amalg\hskip-6pt\amalg(E_{n}/\Q)[2] = 0, $ hence the $
2-$primary part $ \amalg\hskip-6pt\amalg(E_{n}/\Q)[2^{\infty }] = 0,
$ and so $ \sharp \amalg\hskip-6pt\amalg(E_{n}/\Q) $ is odd. But,
from the fact that the full BSD conjecture for $ E_{n} $ over $ \Q,
$ it is easy to know that $ \sharp \amalg\hskip-6pt\amalg(E_{n}/\Q)
= 2^{-2 m} a_{n}^{2} $ (see the above proof of Theorem 3.3), so $
v_{2} (a_{n}) = m. $ The proof is completed. \quad $ \Box $

\par \vskip 0.4 cm

\hspace{-0.8cm}{\bf 4. Application II - Elliptic curves related with
Heegner points.}
\par \vskip 0.2 cm

In this section, let $ E $ be an elliptic curve defined over $ \Q, \
N_{E} $ be the conductor of $ E / \Q, $ let $ K = \Q(\sqrt{D}) $ be
an imaginary quadratic field with fundamental discriminant $ D $
satisfying the Heegner hypothesis, that is,
\par \vskip 0.1 cm
{Heegner hypothesis.} \ All prime numbers $ p $ dividing $ N_{E} $
are split in $ K. $
\par \vskip 0.1 cm
Then there exists a Heegner point $ P_{K} \in E(K) $ (see [GZ],
[Kol1$\sim$3]). We have the following results of Shafarevich-Tate
groups and Heegner points:
\par \vskip 0.2 cm

{\bf Theorem 4.1.} (1) \ Let $ E $ be an elliptic curve defined over
$ \Q, $ and $ K = \Q(\sqrt{D}) $ be an imaginary quadratic field
satisfying the Heegner hypothesis. Let $ P_{K} $ be a Heegner point
of $ E(K), $ if $ P_{K} $ is of infinite order, then
$$ \frac{\sharp \amalg\hskip-6pt\amalg(E/ \Q) \cdot \sharp
\amalg\hskip-6pt\amalg(E_{D}/ \Q)}{\sharp
\amalg\hskip-6pt\amalg(E/K)} = \left \{\begin{array}{l} 2^{1 -
\delta _{\infty } - \delta _{g}} \cdot (E(\Q) : N(E(K)))^{2} \quad
\text{if} \
L(E/ \Q, 1) \neq 0, \\ \\
2^{- 1 - \delta _{\infty } - \delta _{g}} \cdot (E(\Q) :
N(E(K)))^{2} \quad \text{if} \ L(E/ \Q, 1) = 0.
\end{array} \right. $$
(2) \ For the elliptic curve $ E : y^{2} = x^{3} - x + \frac{1}{4} $
and the imaginary quadratic field $ K = \Q(\sqrt{D}) $ satisfying
the Heegner hypothesis, if the Heegner point $ P_{K} \in E(K) $ is
of infinite order, then
$$ \sharp \amalg\hskip-6pt\amalg(E/K) = 2^{\delta _{g}} \cdot \sharp
\amalg\hskip-6pt\amalg(E_{D}/ \Q). $$ In particular, for each $ D
\in \{-7, -11, -47, -71, -83, -84, -127, -159, -164, -219, \\
-231, -263, -271, -287, -292, -303, -308, -359, -371, -404, -443,
-447, -471 \}, $ the group $ \amalg\hskip-6pt\amalg(E/K) $ is
trivial.
\par  \vskip 0.1cm
{\bf Proof.} (1) \ For the elliptic curve $ E / \Q $ and the field $
K, $ by definition, $ S = \{p : \ p \ \text{is a prime number and} \
p \mid D N_{E} \}. $ By the Heegner hypothesis, $ N_{E} $ is prime
to $ D, $ and $ S_{0} = \{p : \ p \ \text{is a prime number and} \ p
\mid D \}, $ in particularly, $ E $ has good reduction at each prime
$ p \in S_{0}, $ so $ S_{g} \cup S_{gu} = S_{0}, $ and then $ S_{a}
= S_{smr} = S_{nsmr} = \emptyset . $ Hence by definition, $ \delta
(E, \Q, K) = \delta _{\infty } + \delta _{g} $ with $ \delta
_{\infty } = 1 $ (resp., 0) if $ \Delta (E) > 0 $ (resp. $ \Delta
(E) < 0 $). On the other hand, by the Heegner hypothesis, from the
functional equation we have $ L(E/K, 1) = 0. $ Since the Heegner
point $ P_{K} $ is of infinite order, by the formula of Gross-Zagier
(see [GZ, p.311]), the analytic rank $ \text{ord}_{s = 1} L(E/K, s)
= 1, $ which implies
\begin{align*} &\text{ord}_{s = 1} L(E/\Q, s) = 1 \quad \text{and}
\quad L(E_{D}/\Q, 1) \neq 0; \quad \text{or} \\
&\text{ord}_{s = 1} L(E_{D}/\Q, s) = 1 \quad \text{and} \quad
L(E/\Q, 1) \neq 0.
\end{align*}
Then by the theorems of Kolyvagin and Gross-Zagier (see
[Kol1$\sim$3], [GZ]), we know that $ r_{\Q} = \text{ord}_{s = 1}
L(E/\Q, s) \quad \text{and} \quad r_{D, \Q} = \text{ord}_{s = 1}
L(E_{D}/\Q, s), $ moreover, all the groups $
\amalg\hskip-6pt\amalg(E/K), \amalg\hskip-6pt\amalg(E/ \Q) $ and $
\amalg\hskip-6pt\amalg(E_{D}/ \Q) $ are finite. The conclusion then
follows from the above refined Yu's formula. This proves (1). \\
(2) \ For the elliptic curve $ E : y^{2} = x^{3} - x + \frac{1}{4},
$ its discriminant $ \Delta (E) = N_{E} = 37 > 0, $ and the equation
$ y^{2} + y = x^{3} - x $ is a global minimal equation of $ E $ over
$ \Q. $ By a theorem of Kolyvagin (see [Kol3, p.444]), we know that
$ L(E/\Q, 1) = 0, \ r_{\Q} = 1 $ and $ \amalg\hskip-6pt\amalg(E/ \Q)
= 0, $ moreover, $ E(\Q) = \Z P_{0} $ with $ P_{0} = (0,
\frac{1}{2}). $ Now from the proof of (1), we have $ \delta _{\infty
} = 1 $ because $ \Delta (E) > 0, $ then by the formula in (1), we
get $$ \sharp \amalg\hskip-6pt\amalg(E/K) = 2^{2 + \delta _{g}}
\cdot (E(\Q) : N(E(K)))^{-2} \cdot \sharp
\amalg\hskip-6pt\amalg(E_{D}/ \Q). \eqno(4.1) $$ Since $ E(\Q) / 2
E(\Q) \cong \Z / 2 \Z \times E(\Q)[2] = \Z / 2 \Z $ because $
E(\Q)[2] = 0, $ by definition, $ (E(\Q) : N(E(K))) \mid (E(\Q) : 2
E(\Q)) = 2, $ hence $ N(E(K)) = E(\Q) $ or $ 2 E(\Q). $ But by the
group law algorithm (see [Si1, p.53]), it is not difficult to verify
that $ P_{0} \notin N(E(K)), $ which implies $ N(E(K)) = 2
E(\Q), $ so $ (E(\Q) : N(E(K))) =2. $ Substituting it into (4.1),
we get
$$ \sharp \amalg\hskip-6pt\amalg(E/K) =
2^{\delta _{g}} \cdot \sharp \amalg\hskip-6pt\amalg(E_{D}/ \Q).
\eqno(4.2) $$ This proves the first conclusion in (2). \\
Now we assume that $ D $ is one of the given 23 integers. Then by a
theorem of Kolyvagin (see [Kol2], p.477), $
\amalg\hskip-6pt\amalg(E_{D}/ \Q) = 0. $ So we only need to compute
$ \delta _{g}. $ From the discussion in (1), we know that $ S_{g}
\cup S_{gu} = S_{0}, $ moreover, it is easy to know that $ S_{gu} =
\{2\} $ if and only if $ D $ is even, otherwise, $ S_{gu} =
\emptyset . $ Furthermore, it can be seen easily that $ E $ has good
supersingular reduction at $ 2. $ Hence by definition, we have
$$ \delta _{g} = \sum _{p \in S_{0} \setminus \{2\}}
\dim _{2} \widetilde{E_{p}}(\F_{p})[2] + \varepsilon (2) $$ with $
\varepsilon (2) = \frac{1}{2}(1 - (-1)^{v_{2}(D)}) $ (resp., $ 0 $)
if $ D $ is even (resp., odd). Obviously, $ \varepsilon (2) = 0 $
for each of these 23 integers, and by calculation, it can be easily
seen that $ \widetilde{E_{p}}(\F_{p})[2] = \{ O \} $ for each $ p
\in S_{0} \setminus \{2\}, $ which implies $ \delta _{g} = 0. $
Therefore, by (4.2), we get \\
$ \sharp \amalg\hskip-6pt\amalg(E/K) = \sharp
\amalg\hskip-6pt\amalg(E_{D}/\Q) = 1, $ that is, $
\amalg\hskip-6pt\amalg(E/K) $ is trivial. This proves (2), and the
proof of is completed. \quad $ \Box $
\par  \vskip 0.2 cm

By the methods of this paper, one can obtain other similar examples
as done in Corollary 3.4 and Theorem 4.1 above.

\par  \vskip 0.3 cm

{ \bf Acknowledgments. } \ I would like to thank the anonymous
referee for a very careful reading of the paper and many helpful
comments and suggestions.

\par  \vskip 0.2 cm

\hspace{-0.8cm} {\bf References }
\begin{description}

\item[[ABF]] J. Antoniadis, M. Bungert, G. Frey, Properties of
twists  of elliptic curves, J. reine angew. math. 405 (1990), 1-28.

\item[[AW]] M. Atiyah, C.T.C. Wall, Cohomology of groups, in:
Algebraic Number Theory (J.W.S. Cassels and A. Frohlich, Eds.),
pp.94-115, London: Academic Press, 1967.

\item[[CW]] J. Coates, A. Wiles, On the conjecture of Birch and
Swinnerton-Dyer, Invent. math., 39 (1977), 223-251.

\item[[GZ]] B.H. Gross, D.B. Zagier, Heegner points and derivatives
of $ L-$series, Invent. math., 84 (1986), 225-320.

\item[[Kob]] N. Koblitz, Introduction to Elliptic Curves and
Modular Forms, 2nd Edition, New York: Springer-Verlag, 1993.

\item[[Kol1]] V.A. Kolyvagin, Finiteness of $ E(\Q) $ and
$ \amalg\hskip-6pt\amalg(E/\Q)$ for a subclass of Weil curves,
(Russian) Izv. Akad. Nauk SSSR Ser. Mat. 52 (1988), 522-540,
670-671; translation in Math. USSR-Izv. 32 (1989), 523-541.

\item[[Kol2]] V.A. Kolyvagin, The Mordell-Weil and Shafarevich-Tate
groups for Weil elliptic curves, (Russian) Izv. Akad. Nauk SSSR Ser.
Mat. 52 (1988), 1154-1180, 1327; translation in Math. USSR-Izv. 33
(1989), 473-499.

\item[[Kol3]] V.A. Kolyvagin, Euler systems. In The Grothendieck
Festschrift, Vol. II, 435-483, Progr. Math. 87, Birkhauser Boston,
Boston, MA, 1990.

\item[[Kr]] K. Kramer, Arithmetic of elliptic curves upon quadratic
extension, Transactions of the American Mathematical Society, 264
(1981), 121-135.

\item[[KT]] K. Kramer, J. Tunnell, Elliptic curves and local $ \varepsilon
-$factors, Compositio Math., 46 (1982), 307-352.

\item[[Ma]] B. Mazur, Rational points of Abelian varieties with
values in towers of number fields, Invent. math., 18 (1972),
183-266.

\item[[R1]] K. Rubin, The main conjectures of Iwasawa theory for
imaginary quadratic fields, Invent. math., 103 (1991), 25-68.

\item[[R2]] K. Rubin, Fudge factors in the Birch and Swinnerton-Dyer
conjecture, in: Ranks of Elliptic Curves and Random Matrix Theory
(J.B. Conrey, D.W. Farmer, F. Mezzadri and N.C. Snaith Eds.),
pp.233-236, Cambridge: Cambridge University Press, 2007.

\item[[RS]] K. Rubin,  A. Silverberg, Rank frequencies for quadratic
twists of elliptic curves, Experiment Math., 10 (2001), 559-569.

\item[[Se]] J. -P. Serre, Local Fields, New York: Springer-Verlag,
1979.

\item[[Si1]] J. H. Silverman, The Arithmetic of Elliptic Curves, GTM
106, 2nd Edition, New York: Springer-Verlag, 2009.

\item[[Si2]] J. H. Silverman, Advanced topics in the Arithmetic
of Elliptic Curves, GTM 151, New York: Springer-Verlag, 1999.

\item[[Ta1]] J. Tate, Duality theorems in Galois cohomology over
number fields. Proc. Intern. Congress Math. at Stockholm, 1962, 288-295.
Institute Mittag-Leffler Djursholm, Sweden, 1963

\item[[Ta2]] J. Tate, Algorithm for determing the type of a singular
fiber in an elliptic pencil, in: Modular functions of one variable,
IV, (Proc. Internat. Summer School, Univ. Antwerp 1972), pp.33-52.
Lecture Notes in Math. 476, Springer, Berlin, 1975.

\item[[T]] J. Tunnell, A classical Diophantine problem and modular
forms of weight $ 3/2, $ Invent. math., 72 (1983), 323-334.

\item[[W]] E. Weiss, Algebraic Number Theory, New York:
McGraw-Hill Book C ompany, Inc, 1963.

\item[[Y]] H. Yu, On Tate-Shafarevich groups over Galois extensions,
Israel J. Math., 141 (2004), 211-220.

\item[[Z]] C. Zhao, A criterion for elliptic curves
with lowest 2-power in L(1), Math. Proc. Cambridge Philos. Soc.
121(1997), 385-400.

\end{description}

\end{document}